\documentclass[onecolumn,notitlepage,superscriptaddress]{revtex4-1}

\usepackage{graphicx}
\usepackage[utf8]{inputenc} 
\usepackage[T1]{fontenc}    
\usepackage{hyperref}       
\usepackage{url}            
\usepackage{booktabs}       
\usepackage{amsfonts}       
\usepackage{nicefrac}       
\usepackage{microtype}      
\usepackage{lipsum}
\usepackage{amsmath,amssymb,amsthm}
\usepackage[ruled]{algorithm2e}
\usepackage{mathtools}
\usepackage{physics}

\usepackage{color}

\newcommand{\calN}{{\mathcal N}}

\newcommand{\bbE}{{\mathbb E}}

\newcommand{\bbR}{{\mathbb R}}

\newcommand{\rmd}{{\mathrm d}}

\newcommand{\rmx}{{\mathrm x}}
\newcommand{\rmy}{{\mathrm y}}
\newcommand{\rmz}{{\mathrm z}}

\newcommand{\tilu}{{\tilde u}}

\newcommand{\tilT}{{\tilde T}}

\newcommand{\tilV}{{\tilde V}}

\makeatletter  
\newcommand{\repeatable}[2]{%
    \global\@namedef{repeatable@#2}{#1}#1 \label{#2}  
}
\newcommand{\repeatref}[1]{%
    \@ifundefined{repeatable@#1}{NOT FOUND}{\footnote[0]{\eqref{#1}$ \displaystyle{ \@nameuse{repeatable@#1} } $}}%
    ~\eqref{#1}} 
\makeatother 

\begin{document}

\title{Model Predictive Mean Field Games for Controlling Multi-Agent Systems}

\author{Daisuke Inoue}
\email{daisuke-inoue@mosk.tytlabs.co.jp}
\affiliation{%
Toyota Central R\&D Labs., Inc.\\
Bunkyo-ku, Tokyo 112-0004, Japan
}%
\author{Yuji Ito}
\affiliation{%
Toyota Central R\&D Labs., Inc.\\
Nagakute, Aichi 480-1192, Japan
}%
\author{Takahito Kashiwabara}
\affiliation{%
Graduate School of Mathematical Sciences, the University of Tokyo\\
3-8-1 Komaba, Meguro-ku, Tokyo 153-8914, Japan
}%
\author{Norikazu Saito}
\affiliation{%
Graduate School of Mathematical Sciences, the University of Tokyo\\
3-8-1 Komaba, Meguro-ku, Tokyo 153-8914, Japan
}%
\author{Hiroaki Yoshida}
\affiliation{%
Toyota Central R\&D Labs., Inc.\\
Bunkyo-ku, Tokyo 112-0004, Japan
}%

\begin{abstract}
  When controlling multi-agent systems, the trade-off between performance and scalability is a major challenge.
  Here, we address this difficulty by using mean field games (MFGs), which is a framework that deduces the macroscopic dynamics describing the density profile of agents from their microscopic dynamics.
  To effectively use the MFG, we propose a model predictive MFG (MP-MFG), which estimates the agent population density profile with using kernel density estimation and manages the input generation with model predictive control. 
  The proposed MP-MFG generates control inputs by monitoring the agent population at each time step, and thus achieves higher robustness than the conventional MFG.
  Numerical results show that the MP-MFG outperforms the MFG when the agent model has modeling errors or the number of agents in the system is small.
\end{abstract}

\maketitle


\section{Introduction}

Multi-agent systems, which consist of many mobile agents, have received attention for their recent progress in robotic and automated driving systems~\cite{Tan2015Handbook,Wooldridge2009introduction}. 
An important aspect of controlling multi-agent systems is to balance achieving control goals with ensuring scalability of the algorithm~\cite{Deloach2001Multiagent,Barca2013Swarm}. 
Here we aim to construct control laws that minimize the evaluation functions associated with each agent 
in a situation where agents have conflicts of interest, that is, their individual profits may be diminished when everyone acts egotistically.
A typical example is a mobility planning problem, in which each vehicle chooses its own route to pursue profit according to the distance to the destination, its own speed, and the degree of congestion.
When everyone acts selfishly and chooses the same path, the congestion increases, resulting in smaller individual profits~\cite{Rapoport2009Choice}.
To avoid such a situation, information on the entire system is necessary to obtain optimal control laws.
This calculation takes an unrealistic amount of time as the number of agents increases~\cite{Huang2020gametheoretic}.

Recently, \emph{mean field games} (MFGs) have been proposed as approximate models for such a many-body optimal control system, and have attracted much attention~\cite{Lasry2007Mean}.
The MFG offers a scalable solution to the aforementioned many-body control problem regardless of the number of agents. 
A key technique in realizing this reduction is to approximate the states of agents as the density field of a group of agents.
This approximation results in two partial differential equations: one equation describing the time evolution of the density distribution of agents, the \emph{Fokker-Planck equation} (FP equation), and the other equation describing the time evolution of the optimal input according to the density evolution, the \emph{Hamilton-Jacobi-Bellman equation} (HJB equation).
Mathematical properties, such as the uniqueness of solutions and conditions for their existence, have been investigated~\cite{bensoussan2013mean,Gueant2011Mean}, and numerical solutions with finite difference methods and finite element methods are available~\cite{achdou2012mean,Lauriere2014Dynamic,Burger2013mean}.
A variety of applications have been achieved using the MFG, such as groups of vehicles on roads~\cite{Huang2020gametheoretic,Festa2018Mean,Tanaka2021Linearly}, groups of pedestrians~\cite{Lachapelle2011mean,Burger2013mean,Dogbe2010Modeling}, swarming robots~\cite{Chen2018Steering}, stock trading~\cite{Cardaliaguet2018Mean}, and infectious diseases~\cite{Djehiche2017MeanFieldType}.
These studies mainly focus on modeling the behavior of optimally controlled population, however, and the implementation of the MFG-based control in real multi-agent systems is still rare.

There are two challenges in actual implementation of agent control using the MFG.
The first one is how to estimate the density profile from the real distribution of the agents.
Because the MFG is a continuous field equation, the estimation of a macroscopic density function from the surrounding agents is necessary. 
In ordinary multi-agent system control, the mean or variance values of the population are used as such macroscopic quantities, but our problem setting is different in that it requires the density function itself.
The second challenge is the performance degradation originating in  the connection between the dynamics of the real agents and that of the density profile.
For example, if modeling errors exist in the real microscopic system, a gap between the predicted population trajectory and the true trajectory occurs, which could result in limited control performance.

\begin{figure}[t]
  \centering
  \includegraphics[width=100mm]{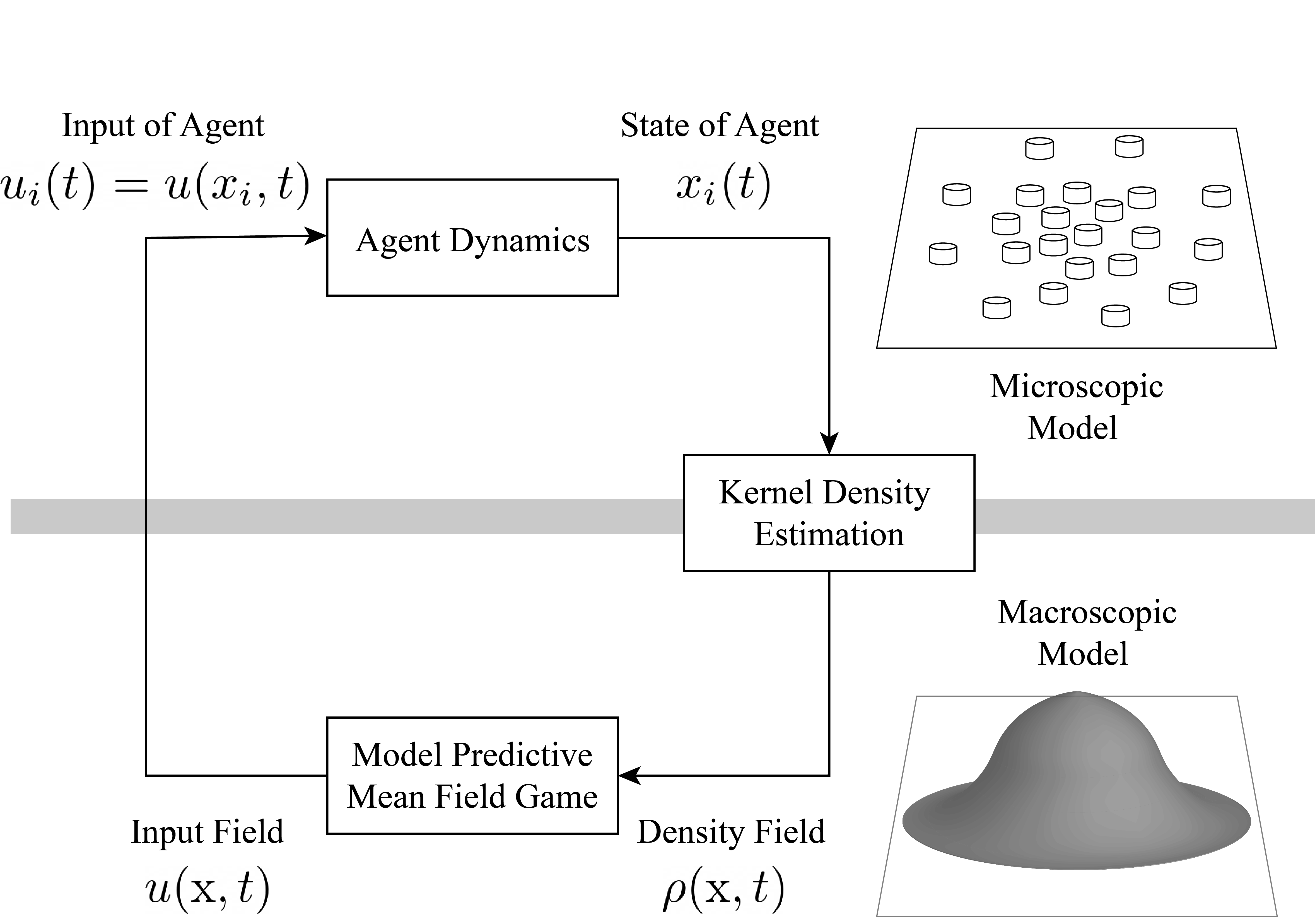}
  \caption{
  Conceptual diagram of the proposed MP-MFG.}
  \label{fig:schematic}
\end{figure}

To overcome these problems, we introduce the following two innovations.
For the first problem, we adopt the \emph{kernel density estimation} (KDE), which is a widely used estimation method.
The KDE has the property that the estimated density approaches the true density as the number of sample points increases~\cite{TurlachBandwidth}.
We propose a method using the KDE in multi-agent systems, with each agent locally performing estimation.
We also give bounds of the estimation error caused by the local observation.
For the second problem, we employ model predictive control.
Specifically, we propose a \emph{model predictive MFG} (MP-MFG) that repeats the operation of 1) estimating the density of the agent population with the KDE and 2) solving the MFG with the estimated density as the initial state at each time step.
Because the prediction error of the density distribution is adjusted by observing the density, the performance degradation is suppressed.
In addition, to reduce the computational complexity, we propose an algorithm that reuses the density information: the MFG solution at previous time step is reused for the initialization procedure of the MFG at the current time step.
A schematic diagram of the MP-MFG is illustrated in Fig.~\ref{fig:schematic}.

The rest of paper is organized as follows.
Section~\ref{sec:problem} formulates the multi-agent control problem, and Section~\ref{sec:MFG} reviews the MFG.
In Section~\ref{sec:proposed}, we propose the MP-MFG control scheme, and in Section~\ref{sec:numeric}, we present numerical experiments to show that the MP-MFG works more robustly than the conventional MFG.
Finally, we provide conclusions and discuss future work in Section~\ref{sec:conclusion}.

\section{Problem Formulation}\label{sec:problem}

Consider $N$ mobile agents located on Euclidean space $\bbR$.
The position of the $i$-th agent at time $t\in [0,T]$, where $T>0$, is denoted as $x_i (t) \in \bbR$, and its dynamics are modeled as an Ito-type stochastic differential equation:
\begin{align}\label{eq:general_dynamics}
  \rmd x_i(t) &= b u_i(t) \rmd t + \sigma \rmd w_i(t),\ t\in[0,T],\\
  x_i(0) &= x_i^0, 
\end{align}
where $b\in\bbR$ and $\sigma\in\bbR$ are the parameters of agents and $x_i^0$ denotes the initial placement of the $i$-th agent. 
The variable $u_i (t) \in \bbR$ is the control input of agent $i$ at time $t$ and $w_i (t) \in \bbR$ is a random variable following the standard Wiener process: $\bbE \{w_i (t) \} = 0 ,\ \bbE \{w_i (t) w_j (s) \} = \delta_{ij} \min(t,s)$.

We define an evaluation function to assess the $i$-th agent input $u_i (t)$ from time $ 0 $ to $T$ as follows:
\begin{align}
  \begin{split}\label{eq:general_eval_func}
  & J_i (\{u_i (t), u_{i^-} (t) \}_{t \in [0, T]} )=\\
  &\quad \bbE\left\{ \int_0^T \left[ \frac{1}{2}(b u_i(t)-\bar v)^2 + q\left(x_i(t), x_{i^-}(t)\right)\right] \rmd t \right\},
  \end{split}
\end{align}
where the first term in the integrand keeps the velocity of the agent close to the desired speed $\bar v\in\bbR$ and the function $q:\bbR\times\bbR^{N-1}\to \bbR$ in the second term evaluates the positional relationship to other agents, which will be defined in more detail later.
The variable $ x_{i ^-} (t) $ is a vector in which the position of all the agents except agent $ i $ are arranged: $x_{i^-}:=[x_1,\ldots, x_{i-1},\ x_{i+1},\ldots, x_{N}]^\top$, and $u_{i ^-} (t)$ is defined in the same fashion.
We attempt to find an input $u_i^*$ that satisfies the following inequality for any $u_i$:
\begin{align}\label{eq:nash_equilibrium}
  J_i(\{u_i^*(t),u_{i^-}^*(t)\}_{t\in[0,T]}) \le  J_i(\{u_i(t),u_{i^-}^*(t)\}_{t\in[0,T]}).
\end{align}
This set of inputs $ u_i^* \ (i = 1, \ldots, N) $ is referred to as a \emph{Nash equilibrium}.
To find the Nash equilibrium, it is necessary to solve an optimal control problem in which Eqs.~\eqref{eq:general_dynamics} and \eqref{eq:general_eval_func} are coupled.
This becomes significantly more difficult as the number of agents $N$ increases.

\section{Review of Mean Field Game}\label{sec:MFG}

Here, we formally derive an approximated model of the aforementioned multi-agent system, that is, the MFG.
A detailed derivation is provided in \cite{Lasry2007Mean, bensoussan2013mean,Gueant2011Mean}.
In the evaluation function of Eq.~\eqref{eq:general_eval_func}, suppose that the function $q$ is rewritten by using a function $\bar q$ as follows:
\begin{align}
  q\left(x_i, x_{i^-}(t)\right) = \bar q(x_i(t),\rho_{i^-}(\cdot,t)), \label{eq:eval_density}
\end{align}
where $\rho_{i^-}$ is the empirical density function defined as
\begin{align}\label{eq:density}
  \rho_{i^-}(\rmx,t) := \frac{1}{N-1}\sum_{j\ne i} \delta(\rmx - x_j(t)),
\end{align}
and $\delta$ denotes the Dirac delta function.
The computational difficulties mentioned above seem to be resolved in Eq.~\eqref{eq:eval_density} because the coupled term $ x_{i ^-} $ has disappeared.
This is incorrect because $ \rho_{i ^-} (\rmx, t) $ actually depends on $ x_{i ^-} $.
However, in the limit where the number of agents is infinitely large, $\rho_{i^-}(\rmx,t)$ is approximated by a continuous density distribution function $\rho(\rmx,t)$ and so is now considered to be a known function.

We consider the following problem, where the second term $q(x_i(t), x_{i^-}(t))$ in Eq.~\eqref{eq:general_eval_func} is replaced by $\bar q(x_i(t), \rho(\rmx,t))$:
\begin{align}
  \begin{split}\label{eq:value_function}
    \min_{u_i([0,T])}\bbE \left\{ \int_t^T \left[ \frac{1}{2}(b u_i(t)-\bar v)^2 + \bar q(x_i(t),\rho(\rmx,t)) \right] \rmd t \right\}.
  \end{split}
\end{align}
Then, it is known that the optimal solution to this problem is obtained by solving the following HJB equation~\cite{Fleming2006Controlled}:
\begin{align}
  \begin{split}
    -\partial_t V(\rmx,t) &= \bar q(\rmx,\rho(\rmx,t)) + \bar v \partial_\rmx V(\rmx,t)\\
    &\quad -\frac{1}{2}(\partial_{\rmx}V(\rmx,t))^2
    + \frac{\sigma^2}{2}\partial_{\rmx\rmx} V(\rmx,t),\label{eq:general_HJB}
  \end{split}\\
  u(\rmx,t) &= \frac{\bar v - \partial_\rmx V (\rmx, t)}{b},\\
  V(\rmx,T) &= 0,\label{eq:general_HJB_terminal}
\end{align}
where $V:\bbR\times[0,T]\to\bbR$ is called a \emph{value function} and $u$ is an optimal input for the problem \eqref{eq:value_function}.
The variables $\partial_{\rmx} {V}(\rmx ,t)$ and $\partial_{\rmx\rmx} {V}(\rmx ,t)$ represent the Jacobian and Hessian, respectively.

Next, the time evolution of the population density is described by the FP equation~\cite{Gardiner2009Stochastic}.
Assume that the control input is given as a known function $ u (\rmx, t) $.
In addition, the initial state of all agents $x_i^0$ is a random variable that follows a predefined probability density function $\rho^0$.
Then, in the limit of $ N \to \infty $, the FP equation representing the time evolution of the density distribution $ \rho(\rmx, t) $ is written as follows:
\begin{align}
  \begin{split}\label{eq:FP}
    \partial_t \rho(\rmx,t) &=  
    - \partial_\rmx \left[ b u(\rmx,t) \rho(\rmx,t)\right] + \frac{\sigma^2}{2}\partial_{\rmx\rmx}\rho(\rmx,t),
  \end{split}\\
  \rho(\rmx,0)&=\rho^0(\rmx).\label{eq:FP_initial_cond}
\end{align}

Finally, a system that combines the HJB equation with the FP equation \eqref{eq:general_HJB}--\eqref{eq:FP_initial_cond} constitutes the MFG.
The MFG includes both the initial condition in the FP equation and the terminal condition in the HJB equation.
We present methods for numerically solving the MFG using the finite element method in Algorithm~\ref{alg:HJB-FP} in the Appendix.

\section{Proposed Method}\label{sec:proposed}

Here, we propose a method to apply the aforementioned MFG to control a multi-agent system.
Our method consists of two ideas: the \emph{KDE}, which estimates the density distribution of the population under an appropriate approximation, and \emph{model predictive control}, which calculates the control input using the MFG at each time step.

\subsection{Kernel Density Estimation}\label{sec:KDE}

For computing the MFG, the density function of the agents is required as an initial condition for the FP equation. 
We adopt the KDE to estimate the density of a group of agents.
Some existing work has used the KDE to control multi-agent systems~\cite{Eren2017Velocitya,Badyn2018Optimal}, but to the best of our knowledge, it has not yet been combined with the MFG.
In the KDE, the estimated density distribution $\hat \rho(\rmx,t)$ is calculated from the agent group location $x_i(t) \ (i=1,\ldots,N)$ as
\begin{align}\label{eq:KDE}
  \hat\rho\left(\rmx, t\right) = \frac{1}{N h} \sum_{j=1}^{N} K\left(\frac{\rmx - x_j(t)}{h}\right),
\end{align}
where $K:\bbR\to\bbR$ is the kernel function and $h>0$ is called the \emph{smoothing parameter}.
The density estimation \eqref{eq:KDE} is known to be asymptotically unbiased, in the sense that the following property holds: 
\begin{align}\label{eq:KDE_converge}
  \lim _{N \rightarrow \infty} \int_{\bbR}\bbE\{\rho(\rmx, t)-\hat\rho(\rmx, t)\}^2\rmd \rmx=0,
\end{align}
when the kernel function $K$ and the smoothing parameter $h$ satisfy suitable conditions~\cite{TurlachBandwidth}.

In order to use the KDE for multi-agent control, we perform a decentralization of the KDE and derive an upper bound on its error.
Assume that each agent $i$ is able to observe the positions of agents $x_j$ within an observable radius $R>0$.
The neighboring agent $\calN_i$ is then defined as
\begin{align}
  \calN_i := \{j\in\{1,\ldots,N\} \mid |x_i(t)-x_j(t)|<R\}.
\end{align}
We evaluate the error between the two estimated density functions, $\check\rho^i(\rmx,t)$ and $\hat\rho(\rmx,t)$, where $\check\rho^i(\rmx,t)$ is the density estimated when agent $i$ performs the KDE using only the locations of neighboring agents and $\hat\rho(\rmx,t)$ is the density estimated using all agents.
We define a local KDE as
\begin{align}
  \check\rho^i\left(\rmx, t\right) = \frac{1}{N h} \sum_{j\in\calN_i} K\left(\frac{\rmx - x_j(t)}{h}\right),
\end{align}
where the kernel function satisfies the following conditions:
\begin{align}
  &K(\rmx)\le K(\rmy)\ \forall |\rmx| > |\rmy|, \label{eq:kernel_prop_1}\\
  & \lim_{|\rmx|\to\infty} K(\rmx)=0, \label{eq:kernel_prop_2}
\end{align}
which are met by many kernel functions~\cite{TurlachBandwidth}.
Then, the estimation error at $\rmx \in \{ \rmz \mid  |\rmz-x_j| \ge R\  \forall j \notin \calN_i \}$
is bounded as follows:
\begin{align}
  \begin{split}\label{eq:kernel_error}
    0 &\le \hat \rho(\rmx,t) - \check\rho^i(\rmx,t)\\
    &= \frac{1}{N h} \sum_{j\notin\calN_i} K\left(\frac{\rmx - x_j(t)}{h}\right)\\
    &\le \frac{1}{N h} \left(N-|\calN_i|\right) K\left(\frac{R}{h}\right),
  \end{split}
\end{align}
where we use Eq.~\eqref{eq:kernel_prop_1}.
When the observable range $R$ is sufficiently large, $|\calN_i|\to N$ holds, and the value on the right side of Eq.~\eqref{eq:kernel_error} approaches 0.
Also, the use of Eq.~\eqref{eq:kernel_prop_2} guarantees that the right-hand side of Eq.~\eqref{eq:kernel_error} becomes sufficiently small when the measurement range $R$ is sufficiently large compared to the smoothing parameter $h$.

\subsection{Model Predictive Mean Field Game}\label{sec:MP-MFG}

When the MFG is used for real system control, modeling errors may accumulate over time and cause deviation between the predicted and actual populations.
To prevent such errors, 
we propose the MP-MFG, in which at each time step $t\in[0,T]$, the system repeatedly solves the MFG of the time interval $ [t, t + \tilT(t)] $ with a prediction horizon $\tilT(t)>0$.
The density function observed at each time step is used as an initial condition, which is expected to suppress the effects of modeling errors.

In the MP-MFG, we discretize the timing of applying the control input with a sampling period $\Delta t>0$ and solve the following MFG with a time window $\tau\in[0,\tilde T(t)]$ for each time step $t\in\{0,\Delta t,\ldots,T\}$:
\begin{align}
  -\partial_\tau \tilV(\rmx,\tau\mid t) &= \bar q(\rmx,\tilde\rho(\rmx,\tau\mid t)) + \bar v \partial_\rmx \tilV(\rmx,\tau\mid t)\label{eq:MP-MFG_HJB}\\
  &\quad -\frac{1}{2}(\partial_{\rmx}\tilV(\rmx,\tau\mid t))^2
  + \frac{\sigma^2}{2}\partial_{\rmx\rmx} \tilV(\rmx,\tau\mid t), \nonumber\\
  \begin{split}
    \partial_\tau \tilde\rho(\rmx,\tau\mid t) &=  
    - \partial_\rmx \left[ b \tilu(\rmx,\tau\mid t) \tilde\rho(\rmx,\tau\mid t)\right]\\
    &\quad + \frac{\sigma^2}{2}\partial_{\rmx\rmx}\tilde\rho(\rmx,\tau\mid t),
  \end{split}\label{eq:MP-MFG_FP}\\
\tilV(\rmx,\tilT(t)\mid t) &= 0,\\
\tilde\rho(\rmx,0\mid t) &= \hat\rho(\rmx,t),\label{eq:MP-MFG_FP_initial}
\end{align}
where $\tilde \rho(\rmx,\tau\mid t)$ and $\tilV(\rmx,\tau\mid t)$ denote the MFG variables at time step $t$ with the time window $ \tau \in [0, \tilT(t)] $.
As shown in Eq.~\eqref{eq:MP-MFG_FP_initial}, the initial density 
$\tilde\rho(\rmx,0\mid t)$
 is set to $\hat\rho(\rmx,t)$, which is estimated with the KDE \eqref{eq:KDE}.
Each agent solves Eqs.~\eqref{eq:MP-MFG_HJB}--\eqref{eq:MP-MFG_FP_initial} at each time step 
$t\in\{0,\Delta t,\ldots,T\}$, and moves with the obtained input:  
\begin{align}\label{eq:MP-MFG_input}
  \tilu (\rmx,0\mid t) = \frac{\bar v -\partial_\rmx \tilV (\rmx,0\mid t)}{b}.
\end{align}

\begin{algorithm}[t]\label{alg:Proposed}
  \caption{Proposed control method: MP-MFG}
  Estimate $\hat\rho(\rmx,0)$ from $x_i^0$ with Eq.~\eqref{eq:KDE}\\
  Calculate $\tilde\rho(\rmx,\tau\mid 0)$ and $\tilV(\rmx,\tau\mid 0)$ for $\tau\in[0,\tilT(0)]$ with Algorithm~\ref{alg:HJB-FP}, using $\tilde\rho(\rmx,\tau\mid 0) = \hat\rho(\rmx,0)$ and $\tilV(\rmx,\tau\mid 0) = 0$ for initial estimation of MFG\\
  $t\leftarrow 0$\\
  \While{$t<T$}{
    Calculate input $\tilu(\rmx,0\mid t)$ with Eq.~\eqref{eq:MP-MFG_input}\\
    Update $x_i(t)$ with $\tilu(\rmx,0\mid t)$\\ 
    $t\leftarrow t+\Delta t$\\
    Construct initial density estimation of MFG $\tilde\rho(\rmx,\tau\mid t)$ for $\tau \in [0,\tilT(t)]$ with Eq.~\eqref{eq:reuse}, using estimated density $\hat\rho(\rmx,t)$ with Eq.~\eqref{eq:KDE} and previous solution of MFG $\tilde\rho(\rmx,\tau\mid t-\Delta t)$\\
    Calculate $\tilde\rho(\rmx,\tau\mid t)$ and $\tilV(\rmx,\tau\mid t)$ for $\tau\in[0,\tilT(t)]$ with Algorithm~\ref{alg:HJB-FP}, using $\tilde\rho(\rmx,\tau\mid t)$ and $\tilV(\rmx,\tau\mid t) = 0$ for initial estimation of MFG
  }
\end{algorithm}

The MP-MFG is computationally expensive because it sequentially solves the MFG, where iterative calculation is performed to find the solution (as shown in Algorithm~\ref{alg:HJB-FP}).
Here, we reduce the computational burden by reusing the density information in the MFG solution between different time steps.
We reuse the predicted density of MP-MFG, $\tilde \rho(\rmx,\tau\mid t-\Delta t)\ (\tau\in[0,\tilde T(t-\Delta t)])$, obtained by finding the input $\tilde u(\rmx,0\mid t-\Delta t)$ at time step $t-\Delta t$.
At the next time step $t$, we initialize the density distribution $\tilde \rho(\rmx,\tau\mid t)\ (\tau\in[0,\tilde T(t)])$ as follows:
\begin{align}
  &\tilde \rho(\rmx,\tau\mid t) =\label{eq:reuse}\\
  \begin{split}    
    &\begin{cases}
      \hat\rho(\rmx, t) & \text{for }\tau = 0,\\
      \tilde \rho(\rmx,\tau+\Delta t\mid t-\Delta t) & \text{for } 0< \tau< \tilT(t-\Delta t)-\Delta t,\\
      \tilde \rho(\rmx,\tilT(t-\Delta t)\mid t-\Delta t) & \text{for } \tau \ge \tilT(t-\Delta t)-\Delta t,
    \end{cases}
  \end{split}\nonumber
\end{align}
where $\hat\rho(\rmx, t)$ is the density estimated at time $t$ using the KDE.
Initializing the density distribution with the predicted density as in Eq.~\eqref{eq:reuse} is expected to speed up the convergence.
A summary of the MP-MFG is shown in Algorithm~\ref{alg:Proposed}.

\section{Numerical Experiment}\label{sec:numeric}

We perform numerical calculations to verify the applicability of the proposed MP-MFG.
The solution is obtained under the domain of periodic interval $\Omega=[0,1]$ and the temporal parameter $T=1$. 
The number of agents is set to $N=1000$, and their initial positions $x_i^0$ are sampled as random variables following a normal distribution with mean $0.2$ and variance $0.1$.
The parameters of the agent dynamics in Eq.~\eqref{eq:general_dynamics} are set to $b=1.0$ and $\sigma=0.1$.
We set the function $\bar q$ in Eq.~\eqref{eq:value_function} as
\begin{align}
  \bar q(\rmx, \rho(\rmx,t)) = C\ln(\rho(\rmx,t)+1),
\end{align}
where $C>0$ is a design parameter.
This function motivates each agent to avoid high-density locations.
The other parameters in the evaluation function are set to $\bar v=1.0$ and $C=0.02$.
In the KDE \eqref{eq:KDE}, the Gaussian kernel~\cite{TurlachBandwidth} is used with the smoothing parameter $h=0.001$.
In the MP-MFG, we use $\tilT(t)=T-t$ for the predictive horizon, $\Delta t=0.001$ for the sampling period, and $\epsilon=3.0\times10^{-6}$ for the parameter to check convergence.

\begin{figure}[t]
  \centering
  \includegraphics[width=100mm]{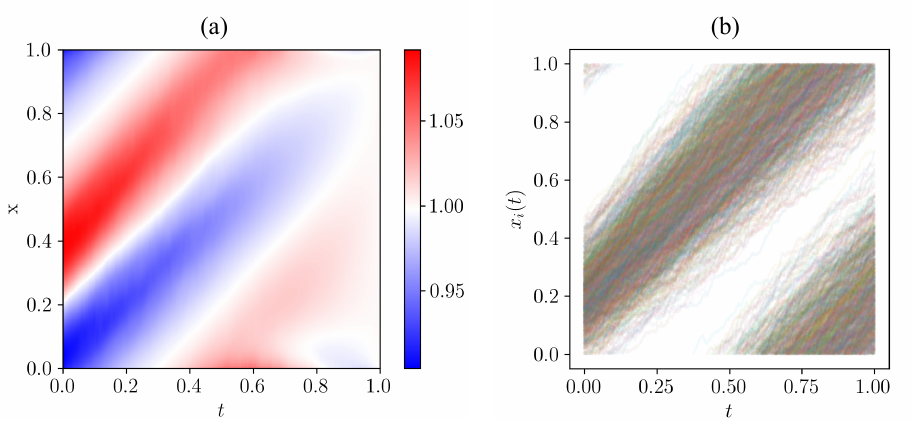}
  \caption{
  Time evolution of the system of Eq.~\eqref{eq:general_dynamics}. 
  (a) Input function $u(\rmx,t)$ was generated with Algorithm~\ref{alg:Proposed}. (b) Position of each agent $x_i(t)$. }
  \label{fig:time-evolution}
\end{figure}

The results of the control with Algorithm~\ref{alg:Proposed} are shown in Fig.~\ref{fig:time-evolution}.
Figure~\ref{fig:time-evolution}(a) shows the obtained input function.
At each time step, inputs that increase or decrease the speed are generated in the region ahead or behind the peak of the density, respectively, indicating that control to avoid concentration is achieved.
Figure~\ref{fig:time-evolution}(b) shows the time response when 1000 agents move with the obtained inputs.
The entire group diffuses because the input causes each agent to avoid high-density regions.
To quantitatively evaluate the performance, we define the average of the evaluation function value as
\begin{align}
  \bar J := \frac{1}{N} \sum_{i=1}^N J_i (\{u_i (t), u_{i^-} (t) \}_{t \in [0, T]} ).
\end{align}
The value of $\bar J$ in the proposed method is $2.08\times 10^{1}$.
In contrast, when all agents behave selfishly, that is, when they use $u_i=\bar v / b$, the value of $\bar J$ is $2.17\times 10^{1}$, indicating that the proposed method provides better performance.

\begin{figure}[t]
  \centering
  \includegraphics[width=100mm]{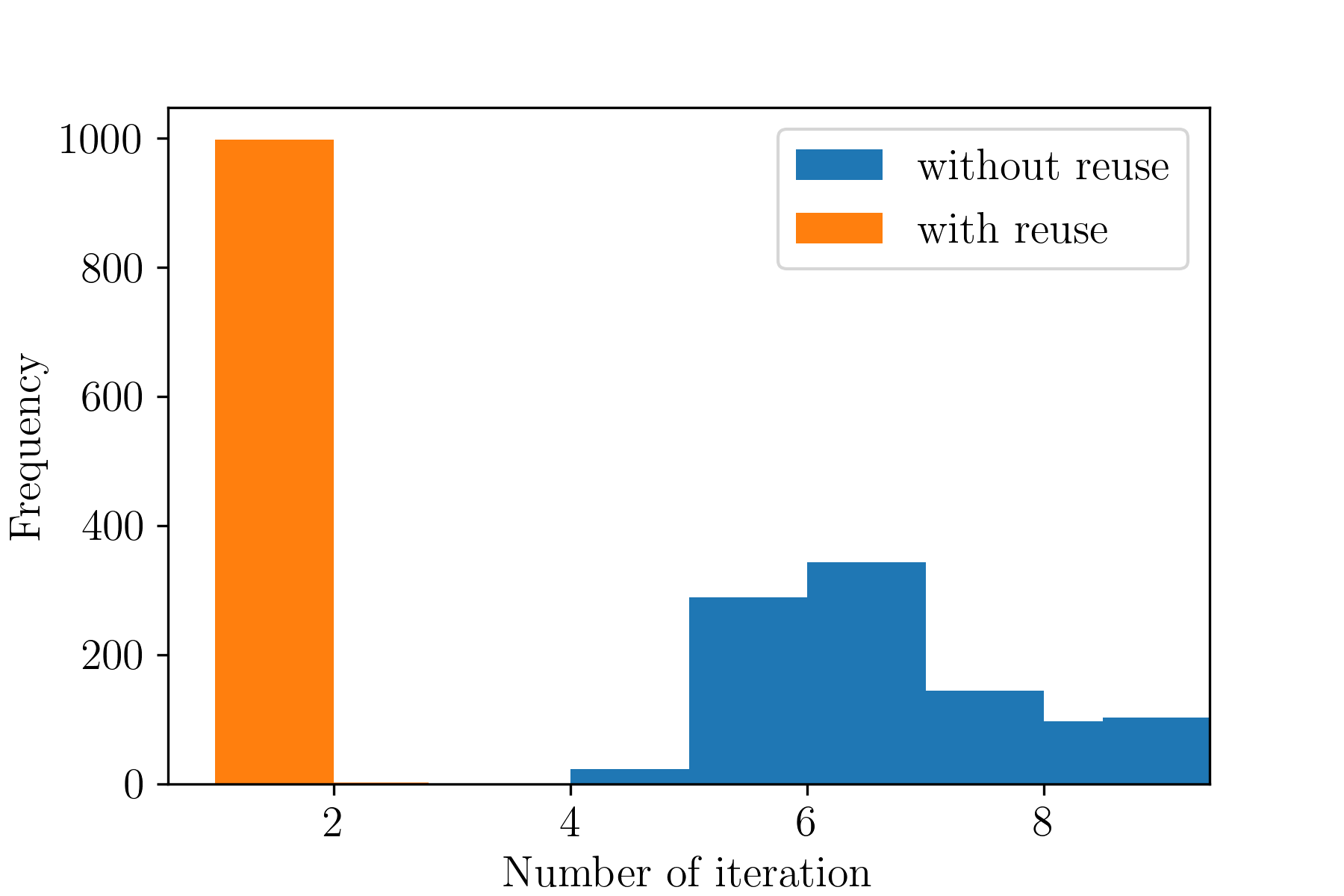}
  \caption{
  Distribution of the number of iterations required for convergence of the MFG.
  The result of reusing the density with Eq.~\eqref{eq:reuse} (orange) is compared with the result of not reusing the density (blue). 
  }
  \label{fig:hist}
\end{figure}

We measure the effect of reusing the density information on the computation time.
In Algorithm~\ref{alg:Proposed}, the number of iterations for running Algorithm~\ref{alg:HJB-FP} at each time is displayed in Fig.~\ref{fig:hist}.
When the density information is not reused and is initialized as $\hat \rho(\rmx,t)$ for all time steps, the average number of iterations is 6.3, while the average value when reusing the density information is 1.0.
Reusing the density information speeds up the calculation by more than 6 times.

\begin{figure}[t]
  \centering
  \includegraphics[width=100mm]{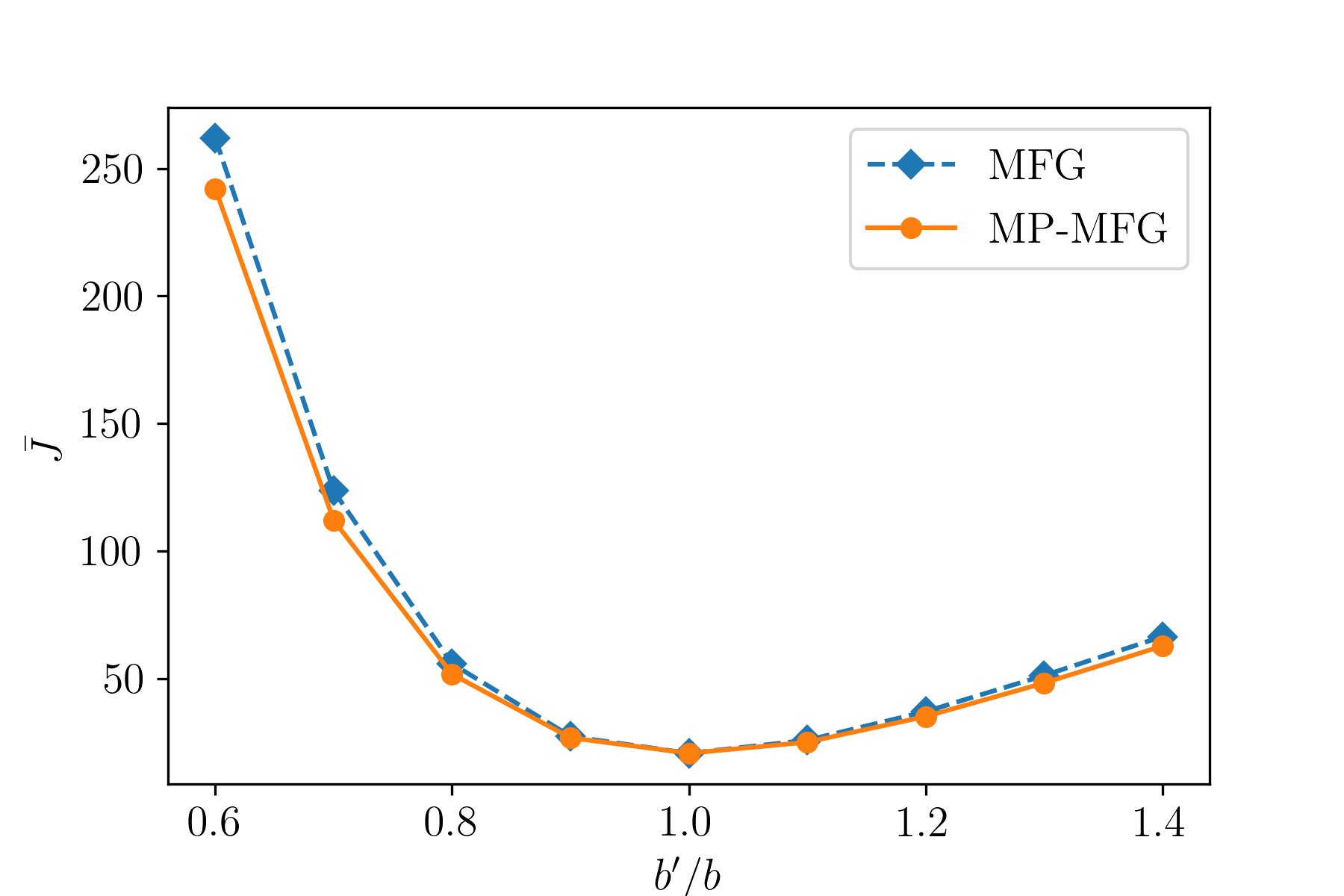}
  \caption{
    Average evaluation function value $\bar J$ for modeling errors $b'/b$.
    The result of the MP-MFG (orange solid curve) is compared with that of the MFG (blue dotted curve). 
  }
  \label{fig:evaluation-function}
\end{figure}

Next, we evaluate the robustness of the proposed method.
When Algorithm~\ref{alg:HJB-FP} is used, the parameter $b'$ is used for the agent dynamics in Eq.~\eqref{eq:general_dynamics} instead of the true parameter $b$.
When the difference between $b$ and $b'$ is large, the control inputs obtained by solving the MFG do not adequately approximate the true solution, and thus the evaluation function value becomes large.
Figure~\ref{fig:evaluation-function} shows the values of $\bar J$ when the modeling error $b'/b$ is changed.
Here, we show both the results of the proposed MP-MFG (Algorithm~\ref{alg:Proposed}) and the conventional MFG (using only Algorithm~\ref{alg:HJB-FP}).
Comparison of the MFG with the MP-MFG shows that the MP-MFG has better performance.
In the MFG, modeling errors cause discrepancies between predicted and actual distributions, resulting in degradation of control performance.
In the proposed MP-MFG, however, the population distribution is estimated by the KDE at each time step, generating relatively reasonable control inputs.

\begin{figure}[t]
  \centering
  \includegraphics[width=100mm]{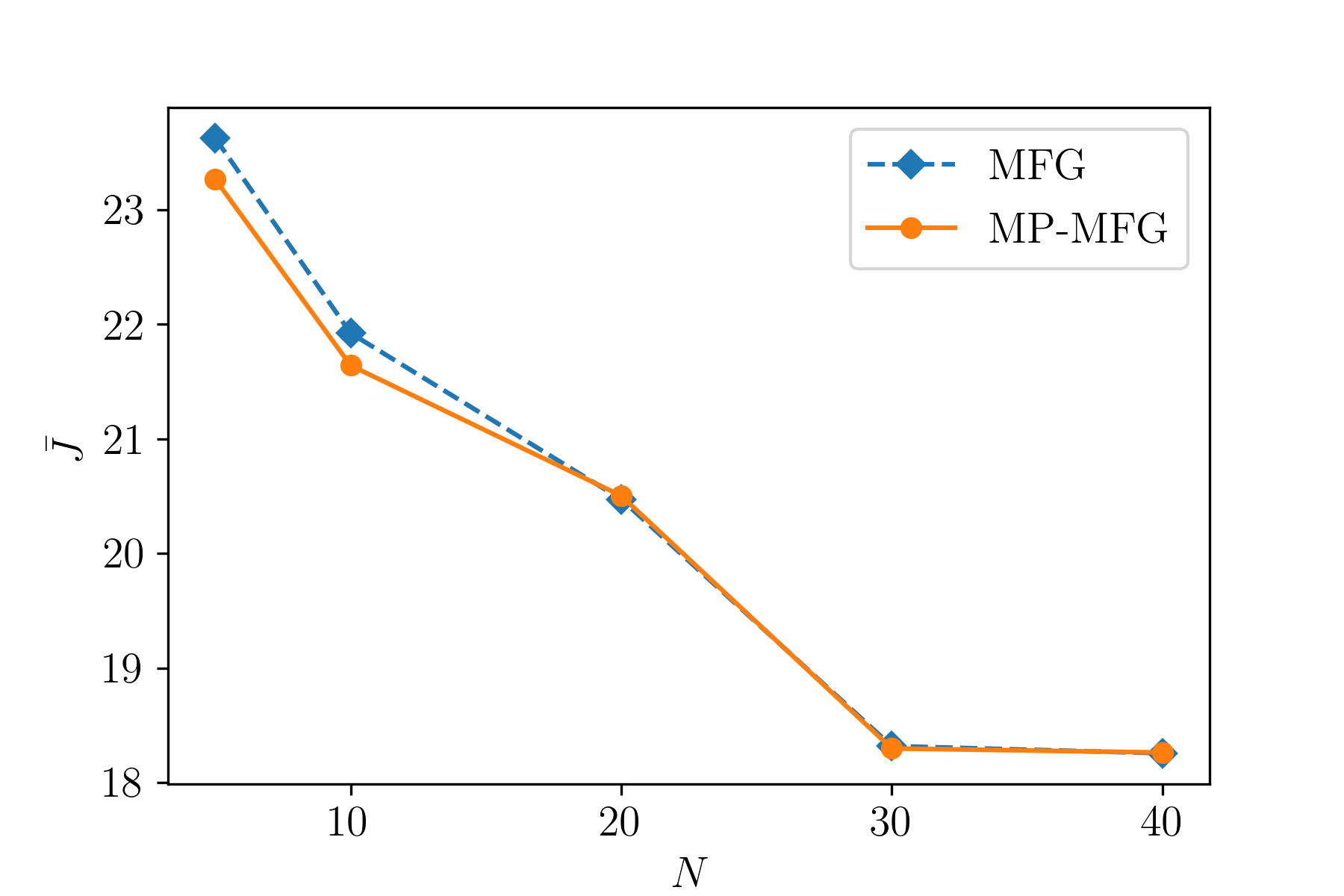}
  \caption{Average evaluation function value $\bar J$ for number of agents $N$.
  The result of the MP-MFG (orange solid curve) is compared with that of the MFG (blue dotted curve). 
  }
  \label{fig:evaluation_num_agent}
\end{figure}

Finally, we evaluate the performance regarding the number of agents.
The values of $\bar J$ for the MFG and MP-MFG for various numbers of agents are shown in Fig.~\ref{fig:evaluation_num_agent}.
When the number of agents is small, the MP-MFG achieves the evaluation function values that are lower than those of the MFG.
The accuracy of the time evolution of a population modeled by the FP equation is degraded due to the small number.
Meanwhile, the MP-MFG re-estimates the density distribution at each time step, so it is possible to track the population relatively accurately, and the performance degradation is suppressed.

\section{Conclusion}\label{sec:conclusion}

In this study, we proposed the MP-MFG for controlling multi-agent systems.
Our method, which combines the KDE and model predictive control, was found to be more robust against modeling errors and better suited for small-scale systems than the control method using conventional MFG.

The effect of the proposed algorithm on the control performance when executed in a distributed manner has not been fully analyzed, and this is a subject for future work.
Also, it is necessary to improve the scheme to counteract the curse of dimensionality when the proposed method is applied to multidimensional systems.

\section*{Appendix}

We derive a numerical scheme for the MFG using the finite element method.
First, the domain of the variable is confined to a interval $\Omega\subset \bbR$ with a proper boundary condition.
Then, the weak forms of Eqs.~\eqref{eq:general_HJB} and \eqref{eq:FP}, which are approximated with time step $\Delta t$, are given as
\begin{align}
  \begin{split}\label{eq:weak_HJB}
    &\int_{\Omega} \frac{V(\rmx ,t+\Delta t) - V(\rmx,t)}{\Delta t}\alpha(\rmx)\rmd \rmx\\
    &+ \int_{\Omega} \bar q(\rmx,\rho(\rmx,t))\alpha(\rmx)\rmd \rmx
    + \int_{\Omega} \bar v \partial_\rmx V(\rmx,t)\alpha(\rmx)\rmd \rmx\\
    &- \int_{\Omega} \frac{1}{2}(\partial_{\rmx}V(\rmx,t))(\partial_{\rmx}V(\rmx,t+\Delta t))\alpha(\rmx)\rmd \rmx\\
    &- \frac{\sigma^2}{2}\int_{\Omega}\partial_{\rmx} V(\rmx,t) \partial_{\rmx}\alpha(\rmx) \rmd \rmx = 0, 
  \end{split}\\
  \begin{split}\label{eq:weak_FP}
    &\int_{\Omega} \frac{\rho(\rmx,t)-\rho(\rmx,t-\Delta t)}{\Delta t} \beta(\rmx)\rmd \rmx\\
    & + \int_{\Omega} bu(\rmx,t) \rho(\rmx,t)\partial_\rmx\beta(\rmx)\rmd \rmx\\
    & + \frac{\sigma^2}{2}\int_{\Omega} \partial_{\rmx}\rho(\rmx,t)\partial_{\rmx}\beta(\rmx)\rmd \rmx= 0, 
  \end{split}
\end{align}
where $\alpha,\beta:\Omega\to\bbR$ are test functions.
After the weak form is obtained, the domain $\Omega$ is divided into finite elements with the width of $\Delta x>0$.
The weak form is then reduced to a finite dimensional linear equation with the values of the physical variables and the value of the test function as the unknown variables.
Finally, these variables are obtained by solving the equation.
These procedures are performed with open-source software FEniCS~\cite{Alnaes2015FEniCS}. 

The calculation of the FP equation can become unstable when the advection term is relatively large compared with the diffusion term.
To avoid this, we add the following stabilization term to the left-hand side of Eq.~\eqref{eq:weak_FP}~\cite{Hughes1989new}:
\begin{align}
  \begin{split}\label{eq:SUPG_FP}
    &\int_{\Omega} \tau \left\{bu(\rmx,t) \partial_\rmx \beta(\rmx) \right\} 
    \gamma(\rmx,t)
    \rmd \rmx,
  \end{split}
\end{align}
where the function $\gamma(\rmx,t)$ is defined as
\begin{align}
  \begin{split}
    \gamma(\rmx,t) &:= \frac{\rho(\rmx,t)-\rho(\rmx,t-\Delta t)}{\Delta t}\\
    & + \partial_\rmx\left[bu(\rmx,t) \rho(\rmx,t)\right] - \frac{\sigma^2}{2} \partial_{\rmx\rmx}\rho(\rmx,t),
  \end{split}
\end{align}
and the variable $\tau$ is defined as
\begin{align}
  \tau = \left\{ \left(\frac{2}{\Delta t}\right)^2 + \left( \frac{2|b u|}{\Delta x} \right)^2 + \left(\frac{2\sigma^2}{\Delta x^2} \right)^2 \right\}^{-\frac{1}{2}}.
\end{align}

\begin{algorithm}[t]\label{alg:HJB}
  \caption{Numerical solution of HJB equation}
  \KwIn{$\rho(\rmx,t)$ and $V(\rmx,T)$}
  \KwOut{$V(\rmx,t)$}
  $t \leftarrow T-\Delta t$, $V(\rmx,t+\Delta t)\leftarrow 0$\\
  \While{$t>0$}{
    Obtain $V(\rmx,t)$ by solving Eq.~\eqref{eq:weak_HJB}, using $V(\rmx, t+\Delta t)$ \\
    $t \leftarrow t - \Delta t$\\
  }
\end{algorithm}

\begin{algorithm}[t]\label{alg:FP}
  \caption{Numerical solution of FP equation}
  \KwIn{$V(\rmx,t)$ and $\rho(\rmx,0)$}
  \KwOut{$\rho(\rmx,t)$}
  $t \leftarrow \Delta t$, $\rho(\rmx,t-\Delta t)\leftarrow \rho(\rmx,0)$\\
  \While{$t<T$}{
    Obtain $\rho(\rmx,t)$ by solving Eqs.~\eqref{eq:weak_FP} and \eqref{eq:SUPG_FP}, using $\rho(\rmx, t-\Delta t)$ \\
    $t \leftarrow t + \Delta t$\\
  }
\end{algorithm}

\begin{algorithm}[t]\label{alg:HJB-FP}
  \caption{Numerical solution of MFG}
  \KwIn{$V_0(\rmx,t)$, $\rho_0(\rmx,t)$: Initial estimation of MFG, $\epsilon>0$: Permissible error}
  \KwOut{$V(\rmx,t)$, $\rho(\rmx,t)$: Solution of MFG}
  $V(\rmx,t)\leftarrow V_0(\rmx,t)$, $\rho(\rmx,t)\leftarrow \rho_0(\rmx,t)$\\
  $z\leftarrow +\infty$\\
  \While{$z>\epsilon$}{
    $V_\text{old}(\rmx,t)\leftarrow V(\rmx,t)$, $\rho_\text{old}(\rmx,t)\leftarrow \rho(\rmx,t)$\\
    Update $V(\rmx,t)$ with Algorithm~\ref{alg:HJB}, using $\rho(\rmx,t)$ and $V(\rmx,T)$ as input\\
    Update $\rho(\rmx,t)$ with Algorithm~\ref{alg:FP}, using $V(\rmx,t)$ and $\rho(\rmx,0)$ as input\\
    $z\leftarrow \|V-V_\text{old}\| + \|\rho-\rho_\text{old}\|$
  }
\end{algorithm}

In Eq.~\eqref{eq:weak_HJB}, the value of $V(\rmx,t)$ at each time step $t$ is calculated backwards using the value of $V(\rmx, t+\Delta t)$. 
Similarly, in Eq.~\eqref{eq:weak_FP}, the value of $\rho(\rmx,t)$ is calculated forward using the value of $\rho(\rmx, t-\Delta t)$. 
The method of updating each equation is described in Algorithms~\ref{alg:HJB} and \ref{alg:FP}. 
We adopt a forward-backward scheme to solve the MFG, where we fix the variables in one equation and update the variables in the other equation. 
This is described in Algorithm~\ref{alg:HJB-FP}.

\end{document}